\documentclass[12pt]{amsart}
\usepackage{amsthm}
\usepackage{epsf}
\usepackage{amssymb}
\usepackage{latexsym}
\usepackage{amsmath}

\newtheorem{theorem}{Theorem}[section] %the resolution could also be [subsectio$

\newtheorem{cor}[theorem]{Corollary}

\newtheorem{prop}[theorem]{Proposition}
\newtheorem{rem}[theorem]{Remark}
\DeclareMathOperator{\Ker}{Ker}

\def\C{{\mathbb C}}
\def\T{{\mathbb T}}
\def\Z{{\mathbb Z}}
\def\P{{\mathbb P}}

\def\semidirect{\ltimes}

\def \ta{\tau}
\def \ta1{\tau_1}

\def \ri{\rightarrow}

\def \CP{\mathbb C \mathbb P}

\newcommand\Gal[1]{{#1}_{\operatorname{Gal}}}
\newcommand\GalAff[1]{{{#1}_{\operatorname{Gal}}^{\operatorname{Aff}}}}

\def\ie{{{i.e.}}}

\def\sub{{\subseteq}}

\newcommand\set[1]{{\{{#1}\}}}
\def\llra{{\;\longrightarrow\;}}
\def\co{{\,{:}\,}}
\def\suchthat{{\,:\,\,}}

\newcommand\figs[1]{Figures/#1}  % Used by Uzi-PC
\def\wrt{{with respect to}}

\def\pitil{{\tilde{\pi}_1}}

\newcommand\sg[1]{{\left<{#1}\right>}}

\def\isom{{\;\cong\;}}
\newcommand\eq[1]{{(\ref{#1})}}
\newcommand\Eq[1]{{Equation (\ref{#1})}}
\newcommand\eqs[2]{{(\ref{#1})--(\ref{#2})}}
\newcommand\Eqs[2]{{Equations (\ref{#1})--(\ref{#2})}}

\newcommand\FIGURE[3][]{{\begin{figure}\epsfysize=0cm {\epsfbox{\figs{#2}}}\caption{#1}\label{#3}\end{figure}}}

\long\def\forget#1\forgotten{}
\long\def\FF#1\FFF{}

\newcommand\begintable[1][] {{}}

\newcommand\CoxY[2][Y]{{\operatorname{C}_{\operatorname{#1}}(#2)}}
\newcommand\ab{{\operatorname{ab}}}

\newcommand\Tref[1]{{Theorem~\ref{#1}}}

\renewcommand\aa[1][r]{{a_{#1}}}
\newcommand\bb[1][r]{{b_{#1}}}
\newcommand\cc[1][r]{{c_{#1}}}
\newcommand\dd[1][r]{{d_{#1}}}
\newcommand\ee[1][r]{{e_{#1}}}
\newcommand\ff[1][r]{{f_{#1}}}
\newif\ifXY % turns XY version on/off
\XYtrue     % Turn it on
\ifXY
\usepackage{xy}
\fi

\ifXY \xyoption{all} \fi

\begin{document}

%% {$$\mbox{Version \VER{}}$$ \vskip-0.18cm}


\date{Received: April 28, 2004; Revised July 2, 2004; Accepted: July 9, 2004}

\title[Coxeter quotient related to $\T \times \T$]
{The Coxeter quotient of the fundamental group of a Galois cover
of $\T \times \T$}

\author[M.~Amram, M.~Teicher, U.~Vishne]
{Meirav Amram$^1$, Mina Teicher and Uzi Vishne}
\stepcounter{footnote} \footnotetext{Partially supported by the
DAAD fellowship (Germany), Eager (Eu-network, HPRN-CT-2009-00099)
and the LDFT postdoctoral fellowship (the Einstein mathematics
institute, Hebrew university, Jerusalem). %} \stepcounter{footnote}
%\footnotetext{
The Emmy Noether Research Institute for
 Mathematics  (center of the Minerva Foundation of Germany), the
Excellency Center "Group  Theoretic Methods in the Study of
Algebraic Varieties of the Israel Science Foundation.}
\stepcounter{footnote}


% \footnotetext{\note{please add here}}

\address{Meirav Amram, the Einstein mathematics institute, Hebrew university, Jerusalem}
\email{ameirav@math.huji.ac.il}
\address{Mina Teicher, department of mathematics, Bar-Ilan university, 52900
Ramat-Gan, Israel} \email{teicher@math.biu.ac.il}
\address{Uzi Vishne, department of mathematics, Bar-Ilan university, 52900
Ramat-Gan, Israel} \email{vishne@math.biu.ac.il}

\begin{abstract}
Let $X$ be the surface $\T\times\T$ where $\T$ is the complex
torus. This paper is the third in a series, studying the
fundamental group of the Galois cover of $X$ \wrt\ a generic
projection onto $\C\P^2$.

Van Kampen Theorem gives a presentation of the fundamental group
of the complement of the branch curve, with $54$ generators and
more than $2000$ relations. Here we introduce a certain natural
quotient (obtained by identifying pairs of generators), prove it
is a quotient of a Coxeter group related to the degeneration of
$X$, and show that this quotient is virtually nilpotent.
\end{abstract}

\maketitle

\section{Overview}\label{overview}

For an algebraic surface $X$ embedded in a projective space
$\C\P^N$, let $\Gal{X}$ be the Galois cover of $X$ with respect to
the full symmetric group. The fundamental group of the Galois
cover is a deformation invariant of surfaces. The first
computation of this invariant can be found in \cite{MoTe1}, where
an algorithm is outlined for the computation of the fundamental
group in terms of generators and relations. Techniques to get a
compact presentation and identify the group are also presented in
this paper, and yet, in the general case it is very difficult to
obtain concrete information on such groups from their
presentation, for example whether the group is virtually solvable.

The group $\pi_1(\Gal{X})$ was computed and identified for
embeddings of the surface $\C\P^1 \times \C\P^1$ \cite{MoTe1}, the
Hirzebruch surfaces (\cite{FRT} and \cite{MRT}), and, recently,
$\C\P^1 \times \T$ \cite{AGTV} where $\T$ is the complex torus.
(See the references therein for more cases).

The computation is done along the following lines. Take a generic
projection $X \ri \C\P^2$ (of degree $n$) with a branch curve $S$.
Let $X_0 \ri \C\P^2$ be the degeneration of $X \ri \C\P^2$ to a
union of planes, and $S_0$ its branch curve in $\C^2$.
The first step is to compute the braid monodromy corresponding to
$S_0$ and use `regeneration rules' \cite{BGT4} to get the braid
monodromy factorization of $S$ (see \cite{AmTe1} for the braid
monodromy notion). Then one applies the van Kampen Theorem
\cite{vk} to get a presentation of $\pi_1(\C^2-S)$ on a standard
set of generators
$\Gamma_1,\Gamma_{1'},\dots,\Gamma_{m},\Gamma_{m'}$, where $2m$ is
the degree of $S$ (see \cite{AmTe2}). % number of lines in $S_0$
Let $\pitil = \pi_1(\C^2-S)/\sg{\Gamma_j^2,\Gamma_{j'}^2}$. There
is a natural
% action of $\pitil$ on the set of planes in the
% degeneration $X_0$ (of size $n$, say), which defines a
homomorphism $\psi \co \pitil \ri S_n$, derived from the natural monodromy
$\pi_1(\C^2 -S) \rightarrow S_n$. %
It is shown in \cite{MoTe1}, that  the kernel of this map is
isomorphic to $\pi_1(\GalAff{X})$, the fundamental group of the
affine part of the Galois cover. Thus we have a short exact
sequence
$$1 \llra \pi_1(\GalAff{X}) \llra \pitil \llra S_n \llra 1.$$
The group $\pi_1(\Gal{X})$ is then obtained by adding the
`projective relation'
\begin{equation}\label{proj}
\Gamma_1 \Gamma_{1'} \dots \Gamma_m \Gamma_{m'} = 1.
\end{equation}

Under the above mentioned monodromy, each pair of generators
$\Gamma_j$ and $\Gamma_{j'}$ is mapped to the same transposition
in $S_n$. Let $C$ denote the quotient of $\pitil$ under the
identification $\Gamma_{j'} = \Gamma_j$ (for all $j$). Taking the
previous sequence modulo the new relation, we get the short exact
sequence
$$1 \llra \pi_1(\GalAff{X})/\sg{\Gamma_j =\Gamma_{j'}} \llra C \llra S_n \llra 1.$$

It is easy to see that $\psi$ splits through $C$, and so we
have the following commutative diagram: %
\ifXY
\begin{equation*}
\xymatrix{%@R=30pt@C=36pt{  %C = length of --- lines.
     1 \ar@{->}[r]
     & \pi_1(\GalAff{X}) \ar@{->}[r] \ar[d]
     & \pitil \ar@{->}[r]^{\psi} \ar[d]
     & S_n \ar@{->}[r] \ar@{=}[d]
     & 1
\\
     1 \ar@{->}[r]
     & \pi_1(\GalAff{X})/\sg{\Gamma_j = \Gamma_{j'}} \ar@{->}[r]
     & C \ar@{->}[r]^{\psi_C}
     & S_n \ar@{->}[r]
     & 1}
\end{equation*}
 \fi

The kernel $K_C = \Ker(\psi_C)$ is then a quotient of
$\pi_1(\Gal{X})$, since the projective relation vanishes when we
identify $\Gamma_j = \Gamma_{j'}$.

For $\C\P^1\times \C\P^1$ and the Hirzebruch surfaces, %
\forget  the degeneration $X_0$ was simply connected. This seems
to be the conceptual explanation for $C$ being isomorphic to $S_n$
in \forgotten %
the group $C$ is isomorphic to $S_n$. \forget (though the details
of this explanation were never written down). \forgotten On the
other hand, for $X = \C\P^1 \times \T$, $C$ was identified to be
the Coxeter group of type $\tilde{A}_5$ (namely isomorphic to $S_6
\semidirect \Z^5$) \cite{AGTV}.

In this paper we obtain a presentation for the group $C$
associated to the surface $X = \T \times \T$, and show that $C$ is
a quotient of a certain Coxeter group (which belongs to the family
that was studied in \cite{v}, which are Coxeter groups with a
natural projection onto a symmetric group, sending the generators
to transpositions). The computation of $\pi_1(\Gal{X})$ for this
surface started in \cite{Am2}, and continued in \cite{AmTe1} and
\cite{AmTe2}. We briefly review in the next section.

Eventually we prove that $K_C$ is abelian % Miriam Beller says: "abelian" is a word, not a name.
by cyclic:
\begin{theorem}\label{main}
% Let $X =\T\times \T$, $S$ the branch curve with respect to a
% certain generic projection to $\C\P^2$
%
The group $C$ is a semidirect product $C = S_{18} \semidirect
K_C$, where $K_C$ is a central extension of $\Z^{34}$ by $\Z$.

More precisely, let $H$ be the group generated by
$x_1,\dots,x_{18}$, $y_{1},\dots,y_{18}$ and $z$, with the
relations
\begin{eqnarray*}
{}[x_i,x_j] & = & 1 , \\
{}[y_i,y_j] & = & 1, \\
{}[x_i,y_j] & = & 1, \\%
{}[x_i,y_i] & = & z
\end{eqnarray*} %
for all $i \neq j$, and $z$ central. Then $K_C$ is isomorphic to
the kernel of the map $H \ri \Z^2 =\sg{x,y}$ defined by $x_i
\mapsto x$ and $y_i \mapsto y$ (and $z \mapsto 1$). The action of
$S_{18}$ on $K_C$ is via indices, and in particular $\sg{z}$ is
the center of $C$.
\end{theorem}

\section{The Coxeter quotient associated to $\T \times \T$}

{}From now on let $X$ denote the surface $\T \times \T$. Let us
recall what has been done in \cite{AmTe1} and \cite{AmTe2}. The
torus $\T$ embeds in $\C\P^2$,  so by the Segre map, $X$ embeds in
$\C\P^{(2+1)(2+1)-1} = \C\P^8$. The torus degenerates as a union
of three lines (in general position), which we depict in %
\ifXY
\begin{figure}[!h]
\begin{equation*}
\xymatrix@R=30pt@C=36pt{  %C = length of --- lines.
     1 \ar@{-}[r]     & 2 \ar@{-}[r]     & 3 \ar@{-}[r]     & 1 }
\end{equation*}
\caption{Degeneration of $\T$}\label{degT}
\end{figure}
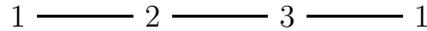%
 \fi
Figure \ref{degT} (with the repeating
index indicating the two points being identified). %
Multiplying two such degenerations, we obtain a degeneration of
$X$ as a union of nine copies of $\C\P^1 \times \C\P^1$, which are
each further degenerated into two copies of $\CP^2$, as seen in
Figure \ref{degenerationTT}. This surface, composed of $18$ planes
with $27$ intersection lines and $9$ intersection points, is
called $X_0$.

\ifXY
\begin{figure}[!h]
\begin{equation*}
\xymatrix@R=30pt@C=36pt{  %C = length of --- lines.
    1 \ar@{--}[r] \ar@{-}[d]|{13}
&
    2 \ar@{--}[r] \ar@{-} [d]|{17} \ar@{-}[ld]|{14}
&
    3 \ar@{--}[r] \ar@{-}[d]|{23} \ar@{-}[ld]|{18}
&
    1 \ar@{--}[d] \ar@{-}[ld]|{22}
\\ %
7 \ar@{-}[r]|{21} \ar@{-}[d]|{15}     & 8 \ar@{-}[r]|{27}
\ar@{-}[d]|{20} \ar@{-}[ld]|{19}     & 9 \ar@{-}[r]|{26}
\ar@{-}[d]| {25} \ar@{-}[ld]|{24}     & 7\ar@{--}[d]
\ar@{-}[ld]|{16} %
\\ %
4 \ar@{-}[r]|{8} \ar@{-}[d]|{4}     & 5 \ar@{-}[r]|{12}
\ar@{-}[d]|{7} \ar@{-}[ld]|{6}     & 6 \ar@{-}[r]|{11}
\ar@{-}[d]|{10} \ar@{-}[ld]|{9}     & 4 \ar@{--}[d]\ar@{-}[ld]|{5}
\\%
     1 \ar@{-}[r]|{1}     & 2 \ar@{-}[r]|{3}     & 3 \ar@{-}[r]|{2}     & 1 }
\end{equation*}
\caption{Degeneration of $X = \T \times \T$}\label{degenerationTT}
\end{figure}%
 \fi

Projecting $X_0$ onto $\C\P^2$, we get a line arrangement $S_0$,
which is the $1$-skeleton of $X_0$, composed of $27$ lines.
Regenerating $X_0$, we get an induced regeneration downstairs from
$S_0$ to the branch curve $S$ of $X$.
In \cite{AmTe1} the degeneration process was described in details,
and the braid monodromy factorization of $S$ was obtained. This
was used in \cite{AmTe2} to compute a presentation
$\pi_1(\C^2-S)$, with the generators
$\Gamma_1,\Gamma_{1'},\dots,\Gamma_{27},\Gamma_{27'}$, and about
$2000$ relations. {}Therefore we have a presentation of $\pitil =
\pi_1(\C^2-S)/\sg{\Gamma_j^2,\Gamma_{j'}^2}$. The generators
correspond (in pairs) to the $27$ lines, and the map $\psi \co
\pitil \ri S_{18}$ is defined by sending $\Gamma_j$ and
$\Gamma_{j'}$ to the transposition $(\alpha \beta)$ where $\alpha$
and $\beta$ are the planes intersecting in line $j$ of $X_0$.

As described above for the general case, the group $C$ is the
quotient of $\pitil$ obtained by adding the relations $\Gamma_{j'}
= \Gamma_j$. We will denote this pair of generators by $u_j$, so
$C =\sg{u_1,\dots,u_{27}}$ with the relations induced from
$\pitil$ under the projection $\theta \co \Gamma_j, \Gamma_{j'}
\mapsto u_j$. Since $\psi(\Gamma_{j'}) = \psi(\Gamma_j)$, $\psi$
splits as the composition $\psi_C \circ \theta$ where $\psi_C$
is defined by %
\begin{equation}%
\label{psiC}\psi_C(u_j) =
\psi(\Gamma_j).%
\end{equation}%
 We remark that the map $u_j \mapsto
\Gamma_j$ does not define a homomorphism from $C$ back to
$\pitil$. The issue of splitting the short exact sequence
$$1 \,\llra\, \Ker(\theta)\,\llra\, \pitil \,\llra\, C \,\llra\, 1$$
is important, but will not be discussed further in this paper.

Our main result was stated as \Tref{main}. For the proof, we first
present $C$ as a quotient of a certain Coxeter group (which is
discussed in the next section), and then apply the general results
of \cite{v}.

\section{Presentations of $S_n$ via transpositions}\label{CYT}

\forget The starting point is the shape of presentations of the
symmetric group on transpositions. We will see that $C$ has most
of the relations of such presentations, and use this observation
to study the kernel of a map $\psi_C \co C \rightarrow S_{18}$.
\forgotten

Let $T$ be an arbitrary graph on $n$ points. To every edge $u \in
T$ we attach the transposition $(\alpha\,\beta)$ where $u$
connects the vertices $\alpha$ and $\beta$. This set of
transpositions generates $S_n$ if and only if $T$ is connected. It
is known that $S_n$ has a presentation with the edges of $T$ as
generators, and the following five sets of relations:
\begin{eqnarray}
{} u^2 & = & 1 \quad \mbox{for all } u \in T, \label{square} \\
{} uv & =& vu \quad \mbox{if $u,v$ are disjoint,} \label{disjoint}\\
{} uvu & = & vuv \quad \mbox{if $u,v$ share a common
vertex,}\label{triple}\\
{} [u,wvw] & = & 1 \quad \mbox{ if $u,v,w$ meet in a common
 vertex,}\label{fork}
\end{eqnarray}
and, for every cycle $u_1,\dots,u_m$ in $T$, the relation
\begin{equation}\label{cycle}
u_1 \dots u_{m-1} = u_2 \dots u_{m},
\end{equation}
which we say is the relation associated to the cycle (see \cite{v}
for details). It is easy to see that (assuming
\eq{square}--\eq{fork} hold) any ordered numeration of the edges
along a cycle gives the same relation \eq{cycle}.
\forget % Misleading, as we are not going to encounter this situation here:
Moreover, it is fairly straightforward that in order to obtain a
presentation for $S_n$, one only needs to assume the relation
\eq{cycle} for a set of generators of the fundamental group
$\pi_1(T)$. \forgotten

Let $\CoxY[]{T}$ denote the Coxeter group generated by $T =
\set{u}$ (one generator for every edge of $T$), with the relations
\eqs{square}{triple}; this is obviously a Coxeter group. Let
$\CoxY{T}$ denote the quotient obtained by adding the relations
\eq{fork}. As we assume $T$ to be connected, the map sending $u$
to the associated transposition is obviously a surjection
\begin{equation}\label{psiT}
\psi_T \co \CoxY{T} \ri
S_n.%
\end{equation}%
We will later show that the group $C$ is a quotient of $\CoxY{T}$
for the graph $T$ of Figure \ref{hexT}. In fact $C$ is obtained by
adding some (but not all) of the cyclic relations \eq{cycle} to
$\CoxY{T}$, so we get a chain of surjections
$$\CoxY[]{T} \llra \CoxY{T} \llra C \llra S_{18}.$$

\section{A presentation of $C$}\label{sec:strucC} % Section 8

The group $C$ was defined above as the image of $\pitil =
\pi_1(\C^2-S)/\sg{\Gamma_j^2,\Gamma_{j'}^2}$ under the map
$\theta$, sending $\Gamma_j$ and $\Gamma_{j'}$ to an abstract
generator $u_j$. A presentation for $\pitil$ was described in
\cite{AmTe2} (with the complete list of relations given in
\cite{Am2}).

Let $T$ denote the graph obtained by connecting the centers of
every two neighboring triangles in Figure~\ref{degenerationTT}.
The resulting graph is given in Figure~\ref{hexT}. Therefore, the
edges of $T$ correspond to lines in $X_0$, and the vertices
correspond to planes. Two edges of $T$ have a joint vertex if and
only if the corresponding lines (in $X_0$) belong to the same
plane (depicted as a triangle in Figure~\ref{degenerationTT}).

\ifXY
\begin{figure}
\begin{equation}\nonumber
\xymatrix@R=5pt@C=6pt{ % C = length of --- lines.
\\
    {}
    &
    &
    &
    &
    &
    &
    &
    &
    &
    &
    &
    &
{}
\\
    {}
    & 1 \ar@{..}[rrr] \ar@{--}[ddd] % \ar@{..}[lllddd]
    &
    &
    & 2 \ar@{..}[rrr] \ar@{--}[ddd] \ar@{--}[lllddd]
    &
    &
    & 3 \ar@{..}[rrr] \ar@{--}[ddd] \ar@{--}[lllddd]
    &
    &
    & 1 \ar@{..}[ddd] \ar@{--}[lllddd]
    &
    &
\\
    {}
    &
    & 7 \ar@{-}[rd] \ar@{-}[ruu] \ar@{-}[lld] %127
    &
    &
    & 9 \ar@{-}[rd] \ar@{-}[ruu] %238
    &
    &
    & 13 \ar@{-}[rd] \ar@{-}[ruu] %139
    &
    &
    &
    &
{}
\\
    {}
    &
    &
    & 11 \ar@{-}[rru] \ar@{-}[ldd]
    &
    &
    & 17 \ar@{-}[rru] \ar@{-}[ldd]
    &
    &
    & 15 \ar@{-}[rru] \ar@{-}[ldd]
    &
    &
    &
{}
\\
    {} %9 \ar@{..}[rrr]
    & 7 \ar@{--}[rrr] \ar@{--}[ddd] %\ar@{..}[lllddd]
    &
    &
    & 8 \ar@{--}[rrr] \ar@{--}[ddd] \ar@{--}[lllddd]
    &
    &
    & 9 \ar@{--}[rrr] \ar@{--}[ddd] \ar@{--}[lllddd]
    &
    &
    & 7\ar@{..}[ddd] \ar@{--}[lllddd]
    &
    &
\\
    {}
    &
    & 12 \ar@{-}[rd] \ar@{-}[lld] %478
    &
    &
    & 18 \ar@{-}[rd] %589
    &
    &
    & 16 \ar@{-}[rd] %679
    &
    &
    &
    &
 {}
\\
    {}
    &
    &
    & 10 \ar@{-}[rru] \ar@{-}[ldd]
    &
    &
    & 14 \ar@{-}[rru] \ar@{-}[ldd]
    &
    &
    & 8 \ar@{-}[rru] \ar@{-}[ldd]
    &
    &
    &
{}
\\
    {} % 6 \ar@{..}[rrr]
    & 4 \ar@{--}[rrr] \ar@{--}[ddd] % \ar@{..}[lllddd]
    &
    &
    & 5 \ar@{--}[rrr] \ar@{--}[ddd] \ar@{--}[lllddd]
    &
    &
    & 6 \ar@{--}[rrr] \ar@{--}[ddd] \ar@{--}[lllddd]
    &
    &
    & 4\ar@{..}[ddd] \ar@{--}[lllddd]
    &
    &
\\
    {}
    &
    & 3 \ar@{-}[rd] \ar@{-}[lld]
    &
    &
    & 6 \ar@{-}[rd]
    &
    &
    & 5 \ar@{-}[rd]
    &
    &
    &
    &
{}
\\
    {}
    &
    &
    & 2 \ar@{-}[rru] \ar@{-}[ldd]
    &
    &
    & 4 \ar@{-}[rru] \ar@{-}[ldd]
    &
    &
    & 1 \ar@{-}[rru] \ar@{-}[ldd]
    &
    &
    &
{}
\\
    {} %    3 \ar@{..}[rrr]
        & 1 \ar@{--}[rrr] % \ar@{..}[ddd]
    &
    &
    & 2 \ar@{--}[rrr] % \ar@{..}[ddd] % \ar@{..}[lllddd]
    &
    &
    & 3 \ar@{--}[rrr] % \ar@{..}[ddd] % \ar@{..}[lllddd]
    &
    &
    & 1 % \ar@{..}[lllddd]
    &
    &
\\
    {}
    &
    &
    &
    &
    &
    &
    &
    &
    &
    &
    &
    &
{} }
\end{equation}
\caption{The graph $T$} \label{hexT}
\end{figure}
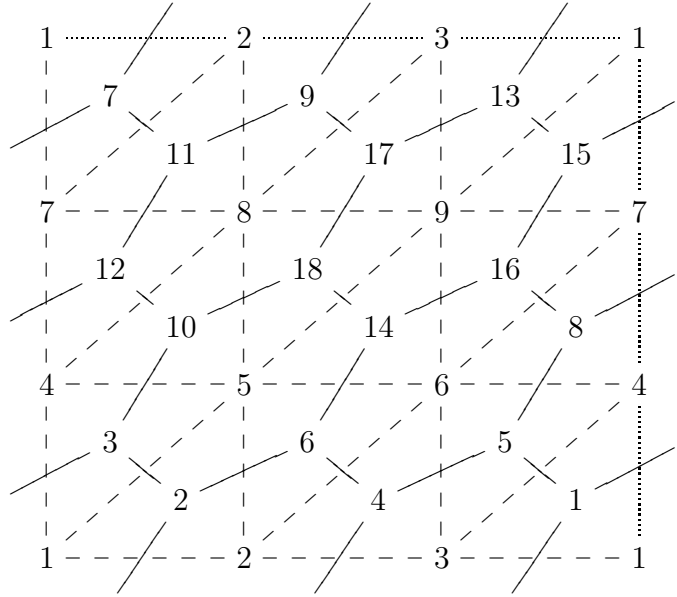
\else \FIGURE{threereg.ps}{three}
\fi %For \ifXY

\begin{theorem}\label{presC}
The group $C$ is generated by $\set{u_j}_{j=1,\dots,27}$, with the
relations \eq{square}--\eq{triple} arising from the graph $T$, and
the cyclic relations \eq{cycle} associated to the nine hexagons in
$T$ (those centered in the intersection points $V_1$--$V_9$ of the
diagram).
\end{theorem}
\begin{proof}
A presentation for $C$ is obtained by substituting $u_j$ for
$\Gamma_j$ and $\Gamma_{j'}$ in the presentation of $\pitil$
(\cite{Am2} and \cite{AmTe2}), which has around 2000 relations.
Fortunately, most of these relations fall into easy to describe
families.

Start with the obvious relations: since $\Gamma_j^2 = 1$ in
$\pitil$, we have $u_j^2 = 1$ in $C$, and so \eq{square} is
proved. Moreover some relations of $C$ have the form $u_i u_j u_i
u_j$, namely $u_i$ commutes with $u_j$. The whole list of
relations was `cleaned' by removing every subword of the form
$u_j^2$, and by replacing every $u_i u_j u_i$ by $u_j$, if $u_i$
and $u_j$ are known to commute. These redundant words were also
removed if they appear in a rotated version of a relation (namely
$rs = 1$ where $sr = 1$ is given). During this process new
commutation relations were `discovered', and they too were used to
further clean the list.
The result of this procedure is a presentation with 333 % was 334, but AX9 was deleted.
relations: 264 commutation relations (of the form $u_iu_j = u_j
u_i$), 44 `triple' relations (of the form $u_i u_j u_i = u_j u_i
u_j$), and
25 %was 26, AX9 deleted.
miscellaneous, which are listed as \Eqs{AX1}{AX26} below (sorted
by length). In order to save space, we write the index $j$ instead
of $u_j$.

\begin{eqnarray}
1 \cdot 13 \cdot 22 \cdot 6 \cdot 4 \cdot 2
\cdot 4 \cdot 6  \cdot 22 \cdot 13 & = & e\label{AX1} \\ % Rotated with respect to original version.
24 \cdot 25 \cdot 26 \cdot 25 \cdot 24 \cdot 22 \cdot 23 \cdot 27 \cdot 23 \cdot 22  & = & e\label{AX2} \\
5 \cdot 4 \cdot 8 \cdot 4 \cdot 5 \cdot 19 \cdot 15 \cdot 11 \cdot
15 \cdot 19  & = & e\label{AX3} \\
9 \cdot 10 \cdot 11 \cdot 10 \cdot 9 \cdot 16 \cdot 25 \cdot 12
\cdot 25 \cdot 16  & = & e\label{AX4}
\end{eqnarray} % length 10
\begin{eqnarray} % length 12
12 \cdot 24 \cdot 20 \cdot 24 \cdot 12 \cdot 24 \cdot 20 \cdot 24 \cdot 12 \cdot 24 \cdot 20 \cdot 24  & = & e\label{AX5} \\
15 \cdot 21 \cdot 15 \cdot 14 \cdot 13 \cdot 14 \cdot 16 \cdot 26
\cdot 16 \cdot 14 \cdot 13 \cdot 14  & = & e\label{AX6} \\
19 \cdot 20 \cdot 19 \cdot 21 \cdot 19 \cdot 20 \cdot 19 \cdot 21 \cdot 19 \cdot 20 \cdot 19 \cdot 21  & = & e\label{AX7} \\
20 \cdot 19 \cdot 21 \cdot 19 \cdot 20 \cdot 18 \cdot 17 \cdot 18 \cdot 27 \cdot 18 \cdot 17 \cdot 18  & = & e\label{AX8} \\
%22 \cdot 23 \cdot 22 \cdot 24 \cdot 25 \cdot 24 \cdot 22 \cdot 23 \cdot 22 \cdot 24 \cdot 25 \cdot 24  & = & e\label{AX9} \\ %-- this is a rotation of \eq{AX10}.
24 \cdot 25 \cdot 24 \cdot 22 \cdot 23 \cdot 22 \cdot 24 \cdot 25 \cdot 24 \cdot 22 \cdot 23 \cdot 22  & = & e\label{AX10} \\
3 \cdot 9 \cdot 7 \cdot 9 \cdot 3 \cdot 9 \cdot 7 \cdot 9 \cdot 3 \cdot 9 \cdot 7 \cdot 9  & = & e\label{AX11} \\
5 \cdot 2 \cdot 5 \cdot 18 \cdot 23 \cdot 18 \cdot 10 \cdot 3
\cdot 10 \cdot 18 \cdot 23 \cdot 18  & = & e\label{AX12} \\
5 \cdot 4 \cdot 5 \cdot 19 \cdot 15 \cdot 19 \cdot 5 \cdot 4 \cdot 5 \cdot 19 \cdot 15 \cdot 19  & = & e\label{AX13} \\
6 \cdot 4 \cdot 6 \cdot 22 \cdot 13 \cdot 22 \cdot 6 \cdot 4 \cdot 6 \cdot 22 \cdot 13 \cdot 22  & = & e\label{AX14}\\
6 \cdot 7 \cdot 6 \cdot 24 \cdot 20 \cdot 24 \cdot 6 \cdot 7 \cdot 6 \cdot 24 \cdot 20 \cdot 24  & = & e\label{AX15} \\
6 \cdot 7 \cdot 6 \cdot 8 \cdot 6 \cdot 7 \cdot 6 \cdot 8 \cdot 6 \cdot 7 \cdot 6 \cdot 8  & = & e\label{AX16} \\
9 \cdot 10 \cdot 9 \cdot 16 \cdot 25 \cdot 16 \cdot 9 \cdot 10 \cdot 9 \cdot 16 \cdot 25 \cdot 16  & = & e\label{AX17} \\
9 \cdot 7 \cdot 9 \cdot 14 \cdot 17 \cdot 14 \cdot 9 \cdot 7 \cdot
9 \cdot 14 \cdot 17 \cdot 14  & = & e\label{AX18}
\end{eqnarray} % length 12
\begin{eqnarray} % length 14
1 \cdot 14 \cdot 17 \cdot 14 \cdot 9 \cdot 7 \cdot 9 \cdot 3 \cdot 9 \cdot 7 \cdot 9 \cdot 14 \cdot 17 \cdot 14  & = & e\label{AX19} \\
6 \cdot 7 \cdot 6 \cdot 8 \cdot 6 \cdot 7 \cdot 6 \cdot 24 \cdot
20 \cdot 24 \cdot 12 \cdot 24 \cdot 20 \cdot 24  & = &
e\label{AX20}
\end{eqnarray}
\begin{eqnarray}
10 \cdot 3 \cdot 10 \cdot 18 \cdot 23 \cdot 18 \cdot 10 \cdot 3
\cdot 10 \cdot 18 \cdot 23 \cdot 18 \cdot 10 \cdot \label{AX21} \\
\cdot 3 \cdot 10
\cdot 18 \cdot 23 \cdot 18  & = & e \nonumber \\
15 \cdot 14 \cdot 13 \cdot 14 \cdot 15 \cdot 21 \cdot 15 \cdot 14
\cdot 13 \cdot 14 \cdot 15 \cdot 21 \cdot 15 \cdot \label{AX22} \\
\cdot 14 \cdot 13 \cdot 14 \cdot 15 \cdot 21  & = & e \nonumber
\end{eqnarray}
\begin{eqnarray}
19 \cdot 20 \cdot 18 \cdot 17 \cdot 18 \cdot 20 \cdot 19 \cdot 21
\cdot 19 \cdot 20 \cdot 18 \cdot 17 \cdot 18 \cdot \label{AX23}\\
\cdot 20 \cdot 19 \cdot 21 \cdot 19 \cdot 20 \cdot 18 \cdot 17
\cdot 18 \cdot 20 \cdot 19 \cdot 21  & = & e \nonumber \\
25 \cdot 24 \cdot 22 \cdot 23 \cdot 22 \cdot 24 \cdot 25 \cdot 26
\cdot 25 \cdot 24 \cdot 22 \cdot 23 \cdot 22 \cdot \label{AX24}\\
\cdot 24 \cdot 25 \cdot 26 \cdot 25 \cdot 24 \cdot 22 \cdot 23
\cdot 22 \cdot 24 \cdot 25 \cdot 26  & = & e \nonumber
\end{eqnarray}
\begin{eqnarray}
6 \cdot 4 \cdot 2 \cdot 4 \cdot 6 \cdot 22 \cdot 13 \cdot 22 \cdot
6 \cdot 4 \cdot 2 \cdot 4 \cdot 6 \cdot 22 \cdot 13 \cdot \label{AX25}\\
\cdot 22 \cdot 6 \cdot 4 \cdot 2 \cdot 4 \cdot 6 \cdot 22 \cdot 13
\cdot 22 & = & e \nonumber \\
9 \cdot 7 \cdot 9 \cdot 3 \cdot 9 \cdot 7 \cdot 9 \cdot 14
\cdot\label{AX26}
17 \cdot 14 \cdot 9 \cdot 7 \cdot 9 \cdot 3 \cdot 9 \cdot \\
7 \cdot 9 \cdot 14 \cdot 17 \cdot 14 \cdot 9 \cdot 7 \cdot 9 \cdot
3 \cdot 9 \cdot 7 \cdot 9 \cdot 14 \cdot 17 \cdot 14  & = &
e\nonumber
\end{eqnarray}

All the 264+44=308 relations of the first two types, which give
the order of products $u_i u_j$ (as $2$ or $3$), match our
expectations: in all the relations $(u_i u_j)^2 = 1$, the edges
$i,j \in T$ do not have a joint vertex, and in all the relations
$(u_i u_j)^3 = 1$, $i$ and $j$ do share a joint vertex. The first
two sets of relations are perhaps best described by listing what
is missing: the $\binom{27}{2} - 308 = 43$ pairs $i,j$ for which
the order is not given in the relations. These $43$
`non-relations' are listed in Table~\ref{nonrel}.

\begin{figure}
\begin{center}
\begin{tabular}{|c|c|c|c|c|}
\hline 1 2  & 1 3  & 1 4  & 1 7 & 1 13     \\
1 17  & 2 3 & 2 10  & 2 13  & 2 23   \\
3 7 & 3 17  & 3 23  & 4 11 & 4 13  \\
4 15  & 7 8 & 7 12 & 7 17 & 7 20  \\
8 11  & 8 12  & 8 15 & 8 20 & 10 12  \\
10 25 & 11 12  & 11 25  & 12 20  & 13 21  \\
13 26  & 15 26  & 17 21 & 17 27  & 20 21 \\
20 27 & 21 26  & 21 27 & 23 25  & 23 26 \\
23 27 & 25 27  & 26 27 & {} & {} \\
\hline
\end{tabular}
\end{center}
\caption{Pairs $i,j$ for which the order of $u_i u_j$ is {\it not
given} as a relation}\label{nonrel}
\end{figure}

It remains to compute the orders of $u_i u_j$ for the pairs of
Table~\ref{nonrel} (thus completing the proof that relations
\eq{disjoint} and \eq{triple} hold), and show that \Eqs{AX1}{AX26}
transform into the nine cyclic relations promised in \eq{cycle}.
Of course these are not all the the cycles in the graph $T$.

Notice that each of the relations \eqs{AX1}{AX26} involves
generators which correspond to lines of $X_0$ with one common
point (namely the edges in $T$ belong to one of the nine
hexagons). Of these, there are nine relations which involve all
the six lines around a point:
\eqs{AX1}{AX4},\eq{AX6},\eq{AX8},\eq{AX12} and \eqs{AX19}{AX20}.
We start by transforming these nine relations into the required
cyclic relations. Some caution is in order here: we consider every
equality of the form $u_m = u_{m-1} \dots u_2 u_1 u_2 \dots
u_{m-1}$ to be a `version' of the cyclic relation; once the orders
of the $u_iu_j$ are known to be $2$ or $3$ (according to whether
or not $i$ and $j$ intersect), all these versions are equivalent.
At this time, however, we do not have all the order relations, and
we will only establish one version of the cyclic relation around
every point.

To simplify reading, we will use the notation given in
Figure~\ref{gen6pt} for the lines around a point: \ifXY
\begin{figure}
\begin{equation*}
%\nonumber
\xymatrix@R=30pt@C=36pt{ % C = length of --- lines.
    {}
    &  \ar@{-}[d]|{\,\bb[]\,}
    &  \ar@{-}[ld]|{\,\cc[]\,}
\\
     \ar@{-}[r]|{\,\aa[]\,}
    & V \ar@{-}[r]|{\,\dd[]\,} \ar@{-}[d]|{\,\ee[]\,} \ar@{-}[ld]|{\,\ff[]\,}
    & {}
\\
    {}
    &  {}
    &  {}
}
\end{equation*}
\caption{Generic names for the lines in $X_0$} \label{gen6pt}
\end{figure}
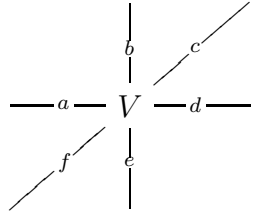
\fi %
If the point is understood from the context, $\aa[]$--$\ff[]$
refer to the appropriate lines around it. In general, we will use
$\aa$--$\ff$ for the lines around the point $V_r$. For example,
$\aa[6] = 12$, $\bb[6] = 25$, and $\dd[5] = 12$ (see
Figure~\ref{degenerationTT}). The advantage of this cumbersome
notation is that we now see that all the missing relations in
Table~\ref{nonrel} involve pairs of the generators $\aa$, $\bb$,
$\dd$ and $\ee$. In fact, in all the pairs $i,j$ of this table,
the lines $i,j$ have a common point; so lines which do no
intersect are known to commute. Likewise diagonals cannot be found
in the table, so every order relation in which $\cc$ or $\ff$ is
involved, is known to hold.

Table~\ref{abde} lists the relations which are not given in
advance.
\begin{figure}
\begin{center}
\begin{tabular}{l|c|c|c|c|c|c|}
      & $\aa\bb$ & $\dd\ee$ & $\bb\dd$ & $\ee\aa$ & $\aa\dd$ & $\bb\ee$ \\
\hline %
$V_1$ &          & x        & x        & x        & x        & x        \\
$V_2$ & x        & x        & x        & x        & x        & x        \\
$V_3$ &          & x        & x        & x        & x        &          \\
$V_4$ &          &          & x        & x        & x        & x        \\
$V_5$ & x        & x        & x        & x        & x        & x        \\
$V_6$ &          &          & x        & x        & x        & x        \\
$V_7$ & x        &          & x        & x        & x        &          \\
$V_8$ & x        & x        & x        & x        & x        &          \\
$V_9$ & x        &          & x        & x        & x        & x        \\
\hline
\end{tabular}
\end{center}
\caption{Order relations which are not given in
advance}\label{abde}
\end{figure}

% In particular:
%   ab de bd ea ad be
% 1     x  x  x  x  x
% 2  x  x  x  x  x  x
% 3     x  x  x  x
% 4        x  x  x  x
% 5  x  x  x  x  x  x
% 6        x  x  x  x
% 7  x     x  x  x
% 8  x  x  x  x  x
% 9  x     x  x  x  x
%
%

Moreover, many relations in \eq{AX1}--\eq{AX26} are instances of
the same `generic' relation in the letters $\aa[]$--$\ff[]$, as we
shall now see.

\begin{rem}\label{trivia}
If $x,y,z$ are elements of order $2$ in a group, and satisfy the
relations $(xy)^3 = (yz)^3 = (xz)^2 = 1$, then $\sg{x,y,z}$ is a
homomorphic image of $S_4$ (which is the Coxeter group of type
$A_3$). In particular $xyzyx = zyxyz$.
\end{rem}

We start with the relations around $V_r$ for $r = 1,4, 6,9$.
Notice that Equations \eq{AX1},\eq{AX3},\eq{AX4} and \eq{AX2} are
the relation $$\ff\ee\dd\ee\ff = \cc\bb\aa\bb\cc$$ around these
four points, respectively. For $r = 4,6,9$, the generators
$\dd,\ee,\ff$ satisfy the relations of Remark~\ref{trivia}, so we
have that $\ff\ee\dd\ee\ff = \dd\ee\ff\ee\dd$, resulting in the
cyclic relation
\begin{equation}\label{cyclic469}
\dd\ee\ff\ee\dd= \cc\bb\aa\bb\cc \qquad \mbox{for $r = 4,6,9$.}
\end{equation}
The situation is slightly different for $r = 1$, since the order
of $\dd[1]\ee[1]$ is not known yet (this is the exception $1,13$
from Table~\ref{nonrel}). However, the orders of
$\aa[1]\bb[1],\aa[1]\cc[1],\bb[1]\cc[1]$ are known, and applying
Remark~\ref{trivia} for $\aa[1],\bb[1],\cc[1]$ we have
\begin{equation}\label{cyclic1}
\ff\ee\dd\ee\ff= \aa\bb\cc\bb\aa \qquad \mbox{for $r = 1$.}
\end{equation}

Next, consider the relation \eq{AX12}, which can be written as
$\cc\dd\cc = \ff\ee\ff\bb\aa\bb\ff\ee\ff$ around $V_3$. Since
$\ff[3] = 18$ commutes with $\cc[3] = 5$ and with $\dd[3] = 2$
(proof: the pairs $2,18$ and $5,18$ cannot be found in
Figure~\ref{nonrel})% This is a perfect proof; we are making the reader used to the nature
% of the group: it is what we define above, nothing more and nothing else.
, we have $\cc[]\dd[]\cc[] = \ee[]\ff[]\bb[]\aa[]\bb[]\ff[]\ee[]$.
Since $\cc[]\dd[]\cc[]=\dd[]\cc[]\dd[]$ (for every $j$) and
$\bb[]\aa[]\bb[]=\aa[]\bb[]\aa[]$ (for $r = 3$), we have that
\begin{equation}\label{cyclic3}
\ee\dd\cc\dd\ee= \ff\aa\bb\aa\ff \qquad \mbox{for $r = 3$.}
\end{equation}
The case of \eq{AX6} is similar: it is the relation
$\ff[]\aa[]\ff[] = \cc[]\bb[]\cc[]\ee[]\dd[]\ee[]\cc[]\bb[]\cc[]$
around $V_7$. Since $\cc[7]$ commutes with $\aa[7]$ and $\ff[7]$,
$\ff[]\aa[]\ff[]=\aa[]\ff[]\aa[]$ (for every $r$), and
$\ee[]\dd[]\ee[]=\dd[]\ee[]\dd[]$ (for $r = 7$), we have that
\begin{equation}\label{cyclic7}
\bb\aa\ff\aa\bb = \cc\dd\ee\dd\cc \qquad \mbox{for $r = 7$.}
\end{equation}

The relations \eq{AX19} and \eq{AX20} both have the form
$$\aa[] = \ff[]\ee[]\ff[]\cc[]\bb[]\cc[]\dd[]\cc[]\bb[]\cc[]\ff[]\ee[]\ff[]$$
around $V_2$ and $V_5$, respectively. We also note that \eq{AX11}
and \eq{AX5} are the relation $\cc[]\bb[]\cc[]\dd[]\cc[]\bb[]\cc[]
= \dd[]\cc[]\bb[]\cc[]\dd[]$ around these points. Combining this
with the fact that $\ff$ commutes with $\bb,\cc$ and $\dd$ (for
every $j$), we have
\begin{equation}\label{cyclic25}
\ee\ff\aa\ff\ee = \dd\cc\bb\cc\dd \qquad \mbox{for $r = 2,5$.}
\end{equation}

Finally, we derive the cyclic relation around $V_8$. Let $x =
\ee[8]\ff[8]\aa[8]\ff[8]\ee[8]$ and $y = \cc[8] \bb[8] \cc[8] =
\bb[8]\cc[8]\bb[8]$. The relation \eq{AX23} is $yxy = xyx$, so
\Eq{AX8}, which is the relation $x = y \dd[8] y$, transforms into
$\dd[] = xyx =
\ee[]\ff[]\aa[]\ff[]\ee[]\bb[]\cc[]\bb[]\ee[]\ff[]\aa[]\ff[]\ee[]$.
But $\ee,\ff$ commute with $\bb,\cc$ (for every $r$, except for
the relation $\bb\ee = \ee\bb$ which does hold for $r = 8$), so we
obtain
\begin{equation}\label{cyclic8}
\ff\ee\dd\ee\ff = \aa\bb\cc\bb\aa \qquad \mbox{for $r = 8$,}
\end{equation}
and this is the last of our nine cyclic relations.

We are now ready to prove relations \eq{square} and \eq{triple}
for the pairs $i,j$ of Figure \ref{nonrel}. As seen from Table
\ref{abde}, the 43 pairs $i,j$ for which the order of $u_iu_j$ is
not given, fall into
four categories: $\aa\dd = \dd\aa$ %two intersecting horizontal lines
(9 pairs, around all the points), $\bb \ee = \ee \bb$ % two intersecting vertical lines
(6 pairs), $\aa\bb\aa=\bb\aa\bb$ and $\dd\ee\dd=\ee\dd\ee$ % horizontal and vertical lines which share a common triangle
(10 pairs%; the corresponding generators should satisfy relation \ref{triple}
), and finally $\bb\dd = \dd \bb$ and $\aa\ee = \ee \aa$ %intersecting horizontal and vertical which do not share a triangle
(18 pairs, two around every point%; the respective generators should commute
). The orders of all the other pairs ($\aa\cc, \aa\ff, \bb\cc,
\bb\ff, \cc\dd, \cc\ee, \cc\ff, \dd\ff$ and $\ee\ff$) are known as
relations for every $r$.

Equations \eq{AX10}, \eq{AX13}, \eq{AX14}, \eq{AX15}, \eq{AX17}
and \eq{AX18} all have the same form, $\cc\bb\cc\ff\ee\ff =
\ff\ee\ff\cc\bb\cc$, around the points $V_r$ with $r = 9,4,1,5,6$
and $2$, respectively. For the other points ($V_3,V_7$ and $V_8$),
we already know that $\bb\ee=\ee\bb$. Since $\ff$ commutes with
$\bb$ and with $\cc$, and $\cc$ commutes with $\ee$, these
relations translate to
\begin{equation}\label{bebe}
\bb\ee=\ee\bb \qquad \mbox{for every $r$.}
\end{equation}

Next, we prove the ten relations of the third kind: that
$u_iu_ju_j = u_ju_iu_j$ if $i,j$ are horizontal and vertical lines
which share a common triangle. The idea is, in each case, to
express $u_i$ (or $u_j$) as a conjugate of another generator using
a cyclic relation, and then show that the conjugate satisfy the
triple relation with $u_j$ (or $u_i$).

As an illustration for this method, consider the pair $\dd[8] =
27$ and $\ee[8] = 20$ (the pair $20,27$ does appear in
Figure~\ref{nonrel}, so the order of $u_{20} u_{27}$ is not yet
known). Notice that $\ee[8] = \bb[5]$. The cyclic relation
\eq{cyclic25} around $V_5$ provides the equality $\bb[5] =
\cc[5]\alpha\cc[5]$ where $\alpha =
\dd[5]\ee[5]\ff[5]\aa[5]\ff[5]\ee[5]\dd[5]$ commutes with $\dd[8]$
(since $\aa[5],\ff[5],\dd[5],\ee[5]$ have no point in common with
$\dd[8]$). As $\dd[8] = \aa[9]$ and $\cc[5] = \ff[9]$, these two
generators satisfy $\dd[8]\cc[5]\dd[8] = \cc[5]\dd[8]\cc[5]$.
Finally, $\dd[8]\ee[8] = \dd[8]\bb[5] =\dd[8]\cc[5]\alpha\cc[5]
\sim \cc[5]\dd[8]\cc[5]\alpha=\dd[8]\cc[5]\dd[8]\alpha \sim % Corrected in Version 3.
\cc[5]\dd[8]\alpha\dd[8]=\cc[5]\alpha= \bb[5]\cc[5]$ (where $\sim$
denotes conjugate in the group), and we are done since
$(\bb[5]\cc[5])^3 = 1$.

The following proposition will be used to prove the relations
$\aa\bb\aa=\bb\aa\bb$, which we need to show for $r = 2,5,7,8,9$.

\begin{prop}\label{prove:ababab}
Let $V_s$ be a point to the left of $V_r$ in $X_0$ (so that
$\aa[r] = \dd[s]$). If $\cc[s]\dd[s]\cc[s]$ commutes with
$\bb[r]$, then $(\aa[r]\bb[r])^3=1$.
\end{prop}
\begin{proof}
Let $\alpha = \cc[s]\dd[s]\cc[s]$. Let $t$ be the index of the
point above $r$, so that $V_r,V_s,V_t$ form a clockwise triangle.
The edges of this triangle are $\aa[r] = \dd[s]$, $\cc[s] =
\ff[t]$ and $\ee[t] = \bb[r]$. Since $\ee[t]\ff[t]\ee[t] =
\ff[t]\ee[t]\ff[t]$, we have $\aa[r]\bb[r] =
 \dd[s]\bb[r] =
 \cc[s]\alpha\cc[s]\bb[r] \sim
 \alpha\cc[s]\bb[r]\cc[s] =
 \alpha\bb[r]\cc[s]\bb[r] \sim
 \bb[r]\alpha\bb[r]\cc[s] =
 \alpha\cc[s] =
 \cc[s]\dd[s]$, and $\cc[s]\dd[s]$ is known to have order $3$.
\end{proof}

The point to the left of $V_r$ for $r = 9$ is $V_s$ for $s = 8$,
and by \eq{cyclic8}, around $V_8$ we have $\dd[] =
\ee[]\ff[]\aa[]\bb[]\cc[]\bb[]\aa[]\ff[]\ee[] =
\ee[]\ff[]\aa[]\cc[]\bb[]\cc[]\aa[]\ff[]\ee[] =
\cc[]\ee[]\ff[]\aa[]\bb[]\aa[]\ff[]\ee[]\cc[]$, so that
$\cc[8]\dd[8]\cc[8]$ is a word in $\aa[8],\bb[8],\ee[8],\ff[8]$
which all commute with $\bb[9]$. By the proposition,
$(\aa[9]\bb[9])^3 = 1$.

For $r = 5,7$ we have $s = 4,9$, respectively. In both cases,
\eq{cyclic469} applies, and we have
$\ff[s]\ee[s]\dd[s]\ee[s]\ff[s] = \dd[s]\ee[s]\ff[s]\ee[s]\dd[s]=
\cc[s]\bb[s]\aa[s]\bb[s]\cc[s]$. But $\cc[s]$ commutes with
$\ee[s],\ff[s]$, so again $\cc[s]\dd[s]\cc[s]$ is a word in the
other letters around $V_s$, which all commute with $\bb[r]$. The
proposition thus gives $(\aa[r]\bb[r])^3 = 1$.

The case $r = 2$ (where $s = 1$) is similar. Around $V_1$ we have
$\dd[] = \ff[]\ee[]\aa[]\bb[]\cc[]\bb[]\aa[]\ee[]\ff[] =
\ff[]\ee[]\aa[]\cc[]\bb[]\cc[]\aa[]\ee[]\ff[]$, but $\cc[]$
commutes with $\aa[],\ee[],\ff[]$, and the same argument applies.

For $r = 8$ the point to the left is $V_s$ for $s = 7$, and the
cyclic relation \eq{cyclic7} is $\bb[7]\aa[7]\ff[7]\aa[7]\bb[7] =
\cc[7]\dd[7]\ee[7]\dd[7]\cc[7]$. But since $\dd[7]\ee[7]\dd[7] =
\ee[7]\dd[7]\ee[7]$ and $\ee[7],\cc[7]$ commute, we find as before
that $\cc[7]\dd[7]\cc[7]$ commutes with $\bb[8]$ and the result
$(\aa[8]\bb[8])^3 = 1$ follows.

Together with the cases $r = 1,3,4,6$ which are given as
relations, we conclude that
\begin{equation}\label{ababab}
\aa\bb\aa=\bb\aa\bb \qquad \mbox{for every $r$.}
\end{equation}

Next, we prove the other half of the third set of relations, \ie\
the relations of the form $\dd\ee\dd=\ee\dd\ee$, which we need to
show for $r = 1,2,3,5,8$. The case $j =8$ was settled above.

\begin{prop}\label{prove:dedede}
Let $V_s$ be a point below $V_r$ (so that $\ee[r] = \bb[s]$). If
$\cc[s]\bb[s]\cc[s]$ commutes with $\dd[r]$, then
$(\dd[r]\ee[r])^3=1$.
\end{prop}
\begin{proof}
As in Proposition \ref{prove:ababab}. Let $\alpha =
\cc[s]\bb[s]\cc[s]$. Let $t$ be the index of the point to the
right of $r$, so that $r,s,t$ form a counterclockwise triangle.
The edges of this triangle are $\ee[r] = \bb[s]$, $\cc[s] =
\ff[t]$ and $\aa[t] = \dd[r]$. Since $\aa[t]\ff[t]\aa[t] =
\ff[t]\aa[t]\ff[t]$, we have $\dd[r]\ee[r] =
 \dd[r]\bb[s] =
 \dd[r]\cc[s]\alpha\cc[s] \sim
 \cc[s]\dd[r]\cc[s]\alpha =
 \dd[r]\cc[s]\dd[r]\alpha \sim
 \cc[s]\dd[r]\alpha\dd[r] = \cc[s]\alpha=
\bb[s]\cc[s]$, and $\bb[s]\cc[s]$ is known to be of order $3$ for
every $s$.
\end{proof}

This proposition immediately applies for $r = 2$ and $r = 5$: in
the first case $s = 8$ and we have from \eq{cyclic8} that
$\cc[s]\bb[s]\cc[s] = \bb[s]\cc[s]\bb[s] =
\aa[s]\ff[s]\ee[s]\dd[s]\ee[s]\ff[s]\aa[s]$, and
$\aa[8],\dd[8],\ee[8],\ff[8]$ all commute with $\dd[2]$. In the
second case $s = 2$ and by \eq{cyclic25} $\cc[s]\bb[s]\cc[s] =
\dd[s]\ee[s]\ff[s]\aa[s]\ff[s]\ee[s]\dd[s]$ and we are done by the
same argument.

The point below $r = 3$ is $s = 9$. The cyclic relation
\eq{cyclic469} provides $\dd[s]\ee[s]\ff[s]\ee[s]\dd[s]=
\cc[s]\bb[s]\aa[s]\bb[s]\cc[s]$. However, using \eq{ababab} we
have $\cc[s]\bb[s]\cc[s] =
\aa[s]\dd[s]\ee[s]\ff[s]\ee[s]\dd[s]\aa[s]$, and the usual
argument applies.

The final case is $r = 1$, where $s = 7$. Then the cyclic relation
\eq{cyclic7} gives $\bb[s]\aa[s]\ff[s]\aa[s]\bb[s] =
\cc[s]\dd[s]\ee[s]\dd[s]\cc[s]$. Again by \eq{ababab} we can apply
Remark~\ref{trivia}, so that $\bb[s]\aa[s]\ff[s]\aa[s]\bb[s] =
\ff[s]\aa[s]\bb[s]\aa[s]\ff[s]$. But $\cc[s]$ commutes with
$\aa[s],\ff[s]$, so $\cc[s]\bb[s]\cc[s] =
\aa[s]\ff[s]\dd[s]\ee[s]\dd[s]\ff[s]\aa[s]$ and Proposition
\ref{prove:dedede} applies (as these generators all commute with
$\dd[r]$). With this, we proved
\begin{equation}\label{dedede}
\dd\ee\dd=\ee\dd\ee \qquad \mbox{for every $r$,}
\end{equation}
and we are done with the third set.

There are three kinds of relations we still need to prove, namely
$(\aa\dd)^2 = 1$, $(\bb\dd)^2 = 1$ and $(\aa\ee)^2=1$, for every
$r$.

In order to prove that $\bb[r],\dd[r]$ commute for every $r$, let
$s$ be the point to the right of $r$, so that $\dd[r] = \aa[s]$.
All we need is to write $\aa[s]$ in terms of the other generators
around $V_s$ (since they have no common point with $\dd[r]$). For
$r = 1,4,5,6,8$, the relations \eq{cyclic469} and \eq{cyclic25}
provide the needed expressions directly. In the other cases the
cyclic relations express $\aa[7]\ff[7]\aa[7]$,
$\aa[3]\bb[3]\aa[3]$ or $\aa[s]\bb[s]\cc[s]\bb[s]\aa[s]$ ($s =
2,8$) in terms of the other generators, but using \eq{ababab} (and
Remark \ref{trivia}) we can write these too as conjugates of the
appropriate $\aa[s]$. Thus we have proved
\begin{equation}\label{bdbdbd}
\bb\dd = \dd\bb \qquad \mbox{for every $r$.}
\end{equation}

In a similar manner we can prove
\begin{equation}\label{aeae}
\aa\ee = \ee\aa \qquad \mbox{for every $r$.}
\end{equation}
Indeed, writing $\aa[r] = \dd[s]$ for an appropriate $s$ ($V_s$ to
the left of $V_r$), all we need is to express $\dd[s]$ in terms of
the other generators around $V_s$. The cyclic relations
\eq{cyclic469}--\eq{cyclic8} express $\dd[s], \dd[s]\cc[s]\dd[s],
\dd[s]\ee[s]\dd[s], \dd[s]\ee[s]\ff[s]\ee[s]\dd[s]$ and
$\dd[s]\cc[s]\bb[s]\cc[s]\dd[s]$ in terms of the other relations
for every $s$; but all these are conjugate to $\dd[s]$ using the
relations we know so far, and we are done.

By now we know the orders of all the products of two generators
around a point, except for $\aa\dd$. In particular, when we
consider $\sg{\bb,\cc,\dd,\ee,\ff}$ or $\sg{\aa,\bb,\cc,\ee,\ff}$
for a fixed $r$, we have a homomorphic image  of the Coxeter group
of type $A_5$, namely the symmetric group $S_6$. The cyclic
relations now present $\aa$ as an element in
$\sg{\bb,\cc,\dd,\ee,\ff}$, which is a transposition disjoint from
$\dd$. For example, around $V_7$ we have (by \eq{cyclic7}) that
$\bb[]\ff[]\aa[]\ff[]\bb[] = \bb[]\aa[]\ff[]\aa[]\bb[] =
\cc[]\dd[]\ee[]\dd[]\cc[] = \ee[]\dd[]\cc[]\dd[]\ee[]$, so that
$\dd[] \aa[] \dd[] =
\dd[]\ff[]\bb[]\ee[]\dd[]\cc[]\dd[]\ee[]\bb[]\ff[]\dd[] =
\ff[]\bb[]\dd[]\ee[]\dd[]\cc[]\dd[]\ee[]\dd[]\bb[]\ff[] =
\ff[]\bb[]\ee[]\dd[]\ee[]\cc[]\ee[]\dd[]\ee[]\bb[]\ff[] =
\ff[]\bb[]\ee[]\dd[]\cc[]\dd[]\ee[]\bb[]\ff[] = \aa[]$. This shows
that
\begin{equation}
\aa\dd = \dd\aa \qquad \mbox{for every $r$,}
\end{equation}
which completes the proof of relations \eq{disjoint} and
\eq{triple}. In particular the remark made after \eq{cycle}
applies, and the cyclic relations we proved become the relations
required in \eq{cycle}.

To finish the proof of the theorem, we need to check that
\Eqs{AX1}{AX26} do not introduce more relations. Since every
relation involves only generators around one point $V_r$, they can
easily be evaluated in $\sg{\aa,\bb,\cc,\dd,\ee,\ff}$, which is a
homomorphic image of the Coxeter group of type $\tilde{A}_5$,
which is isomorphic to $S_6 \semidirect \Z^5$ (in fact the cyclic
relations express $\ff$, say, in terms of the other generators, so
we are computing in homomorphic images of the Coxeter group of
type $A_5$, namely $S_6$). For example, Relation \eq{AX26}
involves generators around $V_2$, and translates to
$(\cc[2]\bb[2]\cc[2]\dd[2]\cc[2]\bb[2]\cc[2] \cdot
\ff[2]\ee[2]\ff[2])^3 = 1$. This can easily be verified in
$\sg{\bb[2],\cc[2],\dd[2],\ee[2],\ff[2]}$ which by now is known to
be a homomorphic image of $S_6$.
\end{proof}

\begin{cor}\label{forktoo}
Relation \eq{fork} is also satisfied by the generators of $C$. In
particular $C$ is a quotient of $\CoxY{T}$.
\end{cor}
\begin{proof}
Suppose that $u_i,u_j,u_k$ are edges of $T$ which meet in a point.
Then $u_j$ and $u_k$ belong to the same hexagon in $T$. Use the
cyclic relation associated to this cycle to rewrite $u_ju_ku_j$ as
a product of generators from the other edges of the cycle, which
in particular to not intersect $u_i$ and therefore commute with
$u_i$.
\end{proof}

\section{The structure of $C$}

\newcommand\lab[1]{x^{(#1)}}

The fundamental group of the graph $T$ is freely generated by $10$
generators. To see this, choose a spanning subtree $T_0$ (which
will contain $18-1 = 17$ edges since $T$ connects $18$ vertices);
then there are $27-17 = 10$ basic cycles, since $T$ has $27$
edges. We label the complement of $T_0$ in $T$ by
$\lab{1},\dots,\lab{10}$, as in Figure~\ref{T0}, where the edges
of the spanning subtree are denoted by double lines. The generator
corresponding to $\lab{\tau}$ is of course the loop resulting from
connecting the end point of $\lab{\tau}$ to the starting point
with the (unique) path on $T_0$.

It is proven in \cite{v} that the natural map from the abstract
group $\CoxY{T_0}$ to $\CoxY{T}$ (sending a generator to itself)
is in fact an embedding. Since $T_0$ is a tree, this group is
isomorphic to the symmetric group on $18$ letters. Moreover, the
cyclic relations defining $C$ as a quotient of $\CoxY{T}$ can be
`solved' in $S_{18}$ (by assigning transpositions to the
generators outside of $T_0$; this will also be evident from the
computations below), and so the subgroup $\sg{u_j \suchthat j \in
T_0}$ of $C$ is isomorphic to $S_{18}$. This constitutes a
splitting of the map $\psi_C \co C \ri S_{18}$.

\ifXY
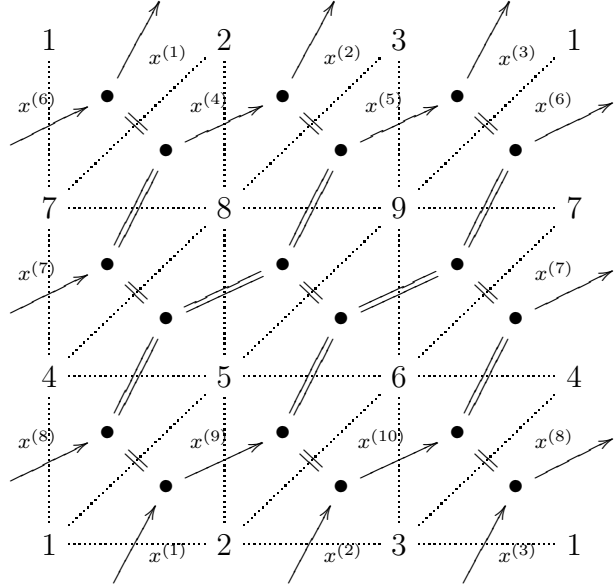
\begin{figure}
\begin{equation}\nonumber
\xymatrix@R=6pt@C=7pt{ % C = length of --- lines.
\\
    {}
    &
    &
    &
    &
    &
    &
    &
    &
    &
    &
    &
    &
{}
\\
    {}
    & 1  \ar@{.}[ddd] % \ar@{..}[lllddd]
    &
    &
    & 2  \ar@{.}[ddd] \ar@{.}[lllddd]
    &
    &
    & 3  \ar@{.}[ddd] \ar@{.}[lllddd]
    &
    &
    & 1  \ar@{.}[lllddd]
    &
    &
\\
    {}
    &
    & \bullet \ar@{=}[rd] \ar@{->}[ruu]_{\lab{1}} \ar@{<-}[lld]_{\lab{6}} %127
    &
    &
    & \bullet \ar@{=}[rd] \ar@{->}[ruu]_{\lab{2}} %238
    &
    &
    & \bullet \ar@{=}[rd] \ar@{->}[ruu]_{\lab{3}} %139
    &
    &
    &
    &
{}
\\
    {}
    &
    &
    & \bullet \ar@{->}[rru]^{\lab{4}} \ar@{=}[ldd]
    &
    &
    & \bullet \ar@{->}[rru]^{\lab{5}} \ar@{=}[ldd]
    &
    &
    & \bullet \ar@{->}[rru]^{\lab{6}} \ar@{=}[ldd]
    &
    &
    &
{}
\\
    {} %9 \ar@{..}[rrr]
    & 7 \ar@{.}[rrr] \ar@{.}[ddd] %\ar@{..}[lllddd]
    &
    &
    & 8 \ar@{.}[rrr] \ar@{.}[ddd] \ar@{.}[lllddd]
    &
    &
    & 9 \ar@{.}[rrr] \ar@{.}[ddd] \ar@{.}[lllddd]
    &
    &
    & 7 \ar@{.}[lllddd]
    &
    &
\\
    {}
    &
    & \bullet \ar@{=}[rd] \ar@{<-}[lld]_{\lab{7}} %478
    &
    &
    & \bullet \ar@{=}[rd] %589
    &
    &
    & \bullet \ar@{=}[rd] %679
    &
    &
    &
    &
 {}
\\
    {}
    &
    &
    & \bullet \ar@{=}[rru] \ar@{=}[ldd]
    &
    &
    & \bullet \ar@{=}[rru] \ar@{=}[ldd]
    &
    &
    & \bullet \ar@{->}[rru]^{\lab{7}} \ar@{=}[ldd]
    &
    &
    &
{}
\\
    {} % 6 \ar@{..}[rrr]
    & 4 \ar@{.}[rrr] \ar@{.}[ddd] % \ar@{..}[lllddd]
    &
    &
    & 5 \ar@{.}[rrr] \ar@{.}[ddd] \ar@{.}[lllddd]
    &
    &
    & 6 \ar@{.}[rrr] \ar@{.}[ddd] \ar@{.}[lllddd]
    &
    &
    & 4 \ar@{.}[lllddd]
    &
    &
\\
    {}
    &
    & \bullet \ar@{=}[rd] \ar@{<-}[lld]_{\lab{8}}
    &
    &
    & \bullet \ar@{=}[rd]
    &
    &
    & \bullet \ar@{=}[rd]
    &
    &
    &
    &
{}
\\
    {}
    &
    &
    & \bullet \ar@{->}[rru]^{\lab{9}} \ar@{<-}[ldd]^{\lab{1}}
    &
    &
    & \bullet \ar@{->}[rru]^{\lab{10}} \ar@{<-}[ldd]^{\lab{2}}
    &
    &
    & \bullet \ar@{->}[rru]^{\lab{8}} \ar@{<-}[ldd]^{\lab{3}}
    &
    &
    &
{}
\\
    {} %    3 \ar@{..}[rrr]
    & 1 \ar@{.}[rrr] % \ar@{..}[ddd]
    &
    &
    & 2 \ar@{.}[rrr] % \ar@{..}[ddd] % \ar@{..}[lllddd]
    &
    &
    & 3 \ar@{.}[rrr] % \ar@{..}[ddd] % \ar@{..}[lllddd]
    &
    &
    & 1 % \ar@{..}[lllddd]
    &
    &
\\
    {}
    &
    &
    &
    &
    &
    &
    &
    &
    &
    &
    &
    &
{} }
\end{equation}
\caption{Spanning subtree of $T$} \label{T0}
\end{figure}
\else \FIGURE{threereg.ps}{three}
\fi %For \ifXY

\forget Let $A_{10,18}$ be the abstract group generated by the
$10\cdot 18^2$ generators $\set{\lab{t}_{ij}}$, $t = 1,\dots,10$,
$i,j =
1,\dots,18$, with the relations %
\begin{eqnarray*}
\lab{t}_{ii} & = & 1\\
\lab{t}_{ij} \lab{t}_{jk} & = & \lab{t}_{ik} \\
\lab{t}_{jk} \lab{t}_{ij} & = & \lab{t}_{ik}
\end{eqnarray*}
for every $\tau = 1,\dots,10$ and $i,j = 1,\dots,18$, and
\begin{equation*}
\lab{t}_{ij} \lab{t'}_{kl} = \lab{t'}_{kl} \lab{t}_{ij}
\end{equation*}
for every $t,t' = 1,\dots,10$ and distinct indices $i,j,k,l \in
\set{1,\dots,18}$. \forgotten

Fix $n = 18$ and $t = 10$. Let $F^{\star}_{t,n}$ be the group
generated by the $t\cdot n$ generators $\set{\lab{\tau}_i}$ ($\tau
= 1,\dots,t$, $i = 1,\dots, n$), with the relations
$$[\lab{\tau}_i, \lab{\tau'}_{j}] = 1 \qquad\mbox{for every $\tau,\tau'$ and
$i \neq j$}$$ %
(therefore $F^{\star}_{t,n}$ is a direct product of $n$ copies of
$\pi_1(T)$ which is the free group on $t$ generators). Let
$e^1,\dots,e^{t}$ denote a set of generators of $\Z^{t}$, and let
$\ab \co F^{\star}_{t,n} \ri \Z^{t}$ be the map defined by
$\ab(\lab{\tau}_i) = e^{\tau}$ (for all $i$). Let $F_{t,n}$ denote
the kernel of this map (note that this is the kernel of the
natural diagonal projection $\pi_1(T)^{n} \ri
\operatorname{H}_1(T)$ since the homology group
$\operatorname{H}_1$ is the abelianization of $\pi_1$).

Recall the definition of $\CoxY{T}$ from Section \ref{CYT}, where
$T$ is the graph of Figure \ref{T0}. Let $N$ denote the normal
subgroup of $\CoxY{T}$ generated by the nine cyclic relations,
associated to the hexagons around the points $V_1,\dots,V_9$. By
Theorem \ref{presC} and Corollary \ref{forktoo}, the group $C$ is
isomorphic to the quotient $\CoxY{T} / N$. The cyclic relations
trivially hold in $S_n$, so the map $\psi_C \co C \ri S_n$ of
\Eq{psiC} is induced from the natural surjection $\psi_T \co
\CoxY{T} \ri S_n$ of \Eq{psiT}.

In \cite{v} (Theorems 5.7 and 6.1) it is shown that $\CoxY{T}
\isom S_{n} \semidirect F_{t,n}$, where $S_{n}$ acts on $F_{t,n}$
by permuting the lower indices. To specify an isomorphism, one
chooses a spanning subtree $T_0$ of $T$ (we take the one given in
Figure \ref{T0}). Then, let $u \in T$ be a (directed) edge,
pointing from $\alpha$ to $\beta$. The isomorphism $\Phi \co
\CoxY{T} \ri S_{n} \semidirect F_{t,n}$ is defined by taking $u$
to the transposition $\Phi(u) = (\alpha\, \beta)$ if $u \in T_0$
(i.e. $u$ is on the spanning subtree), and $\Phi(u) = (\alpha\,
\beta) (\lab{\tau}_{\beta})^{-1}\lab{\tau}_{\alpha}$ if $u =
\lab{\tau}$ is an edge outside of $T_0$. Note that $\Phi(u)$ is an
element of order $2$ in $S_n \semidirect F_{t,n}$. The edges of $T
- T_0$ are ordered for the sake of this definition (since
$(\lab{\tau}_{\beta})^{-1}\lab{\tau}_{\alpha} \neq
(\lab{\tau}_{\alpha})^{-1}\lab{\tau}_{\beta}$ in $F_{t,n}$), but
of course $\lab{\tau}$ and $(\lab{\tau})^{-1}$ is the same element
in $C$ (or even in $\CoxY{T}$).

%%%

Since the cyclic relations hold in the symmetric group, $\psi_T(N)
= 1$ and so $\Phi(N)$ is contained in the kernel of the natural
projection $S_n \semidirect F_{t,n} \ri S_n$, namely $\Phi(N) \sub
F_{t,n}$. Moreover $\Phi(N)$ is normal in $F^{\star}_{t,n}$ so
$$C \isom \CoxY{T}/N \isom (S_n \semidirect F_{t,n})/\Phi(N) = S_n \semidirect (F_{t,n} / \Phi(N))$$
is the kernel of the induced map $\ab \co S_n \semidirect
(F^{\star}_{t,n}/\Phi(N)) \ri \Z^{t}$. We will the quotient
$F^{\star}_{t,n}/\Phi(N)$, and then apply the map $\ab$ and
compute its kernel. See Figure~\ref{groups} for some of the groups
involved. \forget Notice that by definition $\Phi(N)$ is the
normal subgroup generated by (the images under $\Phi$ of) the
cyclic relations, further closed under the action of $S_n$.
\forgotten

\ifXY
\begin{figure}[!h]
\begin{equation*}
\xymatrix@R=30pt@C=36pt{  %C = length of --- lines.
     \CoxY{T} \ar@{->}[r] \ar@{->}[d] \ar@{->}[r]^{\Phi}
     & S_n \semidirect F_{t,n} \ar@{^(->}[r] \ar@{->}[d]
     & S_n \semidirect F^{\star}_{t,n} \ar@{->}[d]
\\
     C \ar@{->}[r] & S_n \semidirect F_{t,n}/\Phi(N) \ar@{^(->}[r]
& S_n \semidirect F^{\star}_{t,n}/\Phi(N)
     }
\end{equation*}
\caption{}\label{groups}
\end{figure}%
 \fi

% proof:
We first compute the image of the cyclic relation associated to
$V_8$. Let $\alpha$ and $\beta$ denote the vertices of $\lab{4}$
(these were planes $11$ and $9$ in Figure~\ref{hexT}). Let
$\lab{4},u_1,\dots,u_5$ denote the edges of the hexagon around
$V_8$ (in that order), then the cyclic relation is that $\lab{4} =
u_1 u_2 u_3 u_4 u_5 u_4 u_3 u_2 u_1$ in $C$ (note that $\lab{4}$
and the $u_j$ have order $2$). Applying $\Phi$, the right hand
side is mapped to the transposition $(\alpha\,\beta)$ while
$\lab{4}$ is mapped to $(\alpha\, \beta)
(\lab{4}_{\beta})^{-1}\lab{4}_{\alpha}$. The equality then becomes
$\lab{4}_{\beta} = \lab{4}_{\alpha}$, which under the action of
$S_n$ becomes $\lab{4}_{j} = \lab{4}_{i}$ for every $i$ and $j$.
Thus $y^4 = \lab{4}_i$ is independent of $i$, and therefore
central (as it commutes with every generator). The same
computation, around $V_5$, $V_6$ and $V_9$, proves that (in $S_n
\semidirect (F^{\star}_{t,n}/\Phi(N))$) $y^5 = \lab{5}_i$, $y^9 =
\lab{9}_i$ and $y^{10} = \lab{10}_i$ are all independent of $i$
and thus central.

Let us now evaluate the cyclic relation around $V_4$. Let
$\lab{7},u_1,u_2,\lab{8}, u_3$ and $u_4$ denote the edges of the
hexagon around $V_4$. Moreover let $\alpha$ and $\beta$ denote the
end points of $\lab{7}$, and $\gamma,\delta$ denote the end points
of $\lab{8}$ ($\lab{8}$ points from $\gamma$ to $\delta$). The
relation in $C$ is
$$\lab{7} u_1 u_2 \lab{8} u_3 = u_1 u_2 \lab{8} u_3 u_4.$$
Applying $\Phi$, we obtain %
$$(\alpha\, \delta) (\lab{7}_{\delta})^{-1} \lab{7}_{\alpha} (\gamma \,
\delta) (\lab{8}_{\delta})^{-1} \lab{8}_{\gamma}  %
= %
(\gamma \,
\delta) (\lab{8}_{\delta})^{-1} \lab{8}_{\gamma}  (\alpha\, \gamma),$$ %
which is equivalent to
$$(\lab{7}_{\delta})^{-1} \lab{7}_{\alpha}   %
= %
 (\lab{8}_{\gamma})^{-1} \lab{8}_{\alpha} (\lab{8}_{\delta})^{-1} \lab{8}_{\gamma},$$
but since $\lab{8}_{\gamma}$, $\lab{8}_{\alpha}$ and
$\lab{8}_{\delta}$ commute, we obtain %
$ \lab{7}_{\alpha}(\lab{8}_{\alpha})^{-1}   %
= %
 \lab{7}_{\delta} (\lab{8}_{\delta})^{-1}$. Acting with $S_n$, we
 obtain
$$\lab{7}_{i}(\lab{8}_{i})^{-1}   %
= %
 \lab{7}_{j} (\lab{8}_{j})^{-1}$$
for every $i,j$. In particular $y^{7,8} =
\lab{7}_{i}(\lab{8}_{i})^{-1}$ is independent of $i$, and
therefore central.

In a similar manner (working around $V_7$, $V_2$ and $V_3$), we
prove that $y^{6,8} = \lab{6}_{i}(\lab{8}_{i})^{-1}$ is
independent of $i$ and central, and likewise for $y^{2,1} =
\lab{2}_{i}(\lab{1}_{i})^{-1}$ and $y^{3,1} =
\lab{3}_{i}(\lab{1}_{i})^{-1}$.

It remains to evaluate the cyclic relation around $V_1$. The
surrounding hexagon is given in Figure~\ref{hex1}, where the
triangles $3,2,7,15,13$ and $1$ (see Figure~\ref{hexT})
were relabelled $\alpha,\beta,\gamma,\delta,\epsilon$ and $\phi$. %

\ifXY
\begin{figure}
\begin{equation}\nonumber
\xymatrix@R=6pt@C=7pt{ % C = length of --- lines.
     &
     &
     &
     & \ar@{.}[ddd]
     &
     &
     & \ar@{.}[lllddd]
\\
    {}
    &
    &
    &
    &
    & \alpha \ar@{=>}[rd]^{u_1} %589
    &
    &
    &
    &
    &
 {}
\\
    {}
    &
    &
    & \phi \ar@{->}[rru]^{\lab{8}} \ar@{<-}[ldd]_{\lab{3}}
    &
    &
    & \beta \ar@{<-}[ldd]^{\lab{1}}
    &
    &
    &
    &
{}
\\
    {} % 6 \ar@{..}[rrr]
    &  \ar@{.}[rrr]  % \ar@{..}[lllddd]
    &
    &
    & V_1 \ar@{.}[rrr] \ar@{.}[ddd] \ar@{.}[lllddd]
    &
    &
    &
    &
    &
    &
\\
    {}
    &
    & \epsilon \ar@{<=}[rd]^{u_2}
    &
    &
    & \gamma
    &
    &
    &
    &
    &
{}
\\
    {}
    &
    &
    & \delta \ar@{->}[rru]_{\lab{6}}
    &
    &
    &
    &
    &
    &
    &
{} \\
     & & & & & & &
 }
\end{equation}
\caption{The hexagon around $V_1$} \label{hex1}
\end{figure}
\else
\fi %For \ifXY

The cyclic relation around $V_1$ in $C$ is %
$$u_1 \lab{1} \lab{6} u_2 \lab{3} = \lab{1} \lab{6} u_2
\lab{3} \lab{8}.$$ %
Applying $\Phi$, we obtain %
\begin{eqnarray*} & & (\alpha\,\beta)
(\beta\,\gamma) (\lab{1}_{\beta})^{-1} \lab{1}_{\gamma}
(\gamma\,\delta) (\lab{6}_{\gamma})^{-1} \lab{6}_{\delta}
 (\delta\,\epsilon)
 (\epsilon\,\phi) (\lab{3}_{\phi})^{-1} \lab{3}_{\epsilon} \\
& = & %
(\beta\,\gamma) (\lab{1}_{\beta})^{-1} \lab{1}_{\gamma}
(\gamma\,\delta) (\lab{6}_{\gamma})^{-1} \lab{6}_{\delta}
(\delta\,\epsilon)
 (\epsilon\,\phi) (\lab{3}_{\phi})^{-1} \lab{3}_{\epsilon}
 (\phi\,\alpha) (\lab{8}_{\alpha})^{-1} \lab{8}_{\phi},\end{eqnarray*}%
which translates to
\begin{equation*}
(\lab{1}_{\gamma})^{-1} \lab{1}_{\beta} (\lab{6}_{\delta})^{-1}
\lab{6}_{\beta}
 (\lab{3}_{\beta})^{-1} \lab{3}_{\phi}
= (\lab{1}_{\gamma})^{-1} \lab{1}_{\alpha} (\lab{6}_{\delta})^{-1}
\lab{6}_{\alpha} (\lab{3}_{\alpha})^{-1} \lab{3}_{\phi}
  (\lab{8}_{\alpha})^{-1} \lab{8}_{\beta},\end{equation*}%
and then (using commutation) to
\begin{equation*}
 \lab{1}_{\beta}
\lab{6}_{\beta}
 (\lab{3}_{\beta})^{-1} (\lab{8}_{\beta})^{-1}
 =
 \lab{1}_{\alpha}  \lab{6}_{\alpha} (\lab{3}_{\alpha})^{-1}
 (\lab{8}_{\alpha})^{-1} .\end{equation*}%
But $\lab{6}_i (\lab{8}_i)^{-1} = y^{6,8}$ and $\lab{3}_i
(\lab{1}_i)^{-1} = y^{3,1}$, are central, so acting with $S_n$, we
obtain
\begin{equation*}
 \lab{1}_{j}
 \lab{8}_{j}
(\lab{1}_{j})^{-1} (\lab{8}_{j})^{-1}  = \lab{1}_{i} \lab{8}_{i}
(\lab{1}_{i})^{-1}
 (\lab{8}_{i})^{-1}
  \end{equation*}%
for every $i,j$. It follows that $z = [\lab{1}_{i}, \lab{8}_{i}]$
is independent of $i$ (and therefore central). %

Summarizing, $F^{\star}_{t,n} / \Phi(N)$ is generated by $y^4,
y^5, y^9, y^{10}$, $y^{7,8}, y^{6,8}, y^{3,1}, y^{2,1}$ and $z$
which are all central, and by $\set{\lab{1}_i, \lab{8}_i}_{i =
1,\dots,18}$, subject to the relations
\begin{eqnarray}
{}[\lab{1}_i,\lab{1}_j] & = & 1 \label{Fs0}, \\
{}[\lab{8}_i,\lab{8}_j] & = & 1, \\
{}[\lab{1}_i,\lab{8}_j] & = & 1, %
\end{eqnarray}
(for all $i \neq j$)  and
\begin{equation} %
{}[\lab{1}_i,\lab{8}_i] = z \label{Fs1},
\end{equation} %
for all $i$. Chasing back the definition of the various $y$
generators, we see that the map $\ab \co F^{\star}_{t,n} \ri
\Z^t=\set{e^1,\dots,e^t}$ is defined by $\ab(y^{\tau}) = e^{\tau}$
for $\tau = 4,5,9,10$, $\ab(y^{\tau,\tau'}) =
e^{\tau}(e^{\tau'})^{-1}$ for $(\tau,\tau') =
(6,8),(7,8),(3,1),(2,1)$, $\ab(\lab{\tau}_i) = e^{\tau}$ for $\tau
= 1,8$ and $\ab(z) = 1$. Let $g \in F^{\star}_{t,n}$ be an
arbitrary element. For every $\tau\neq 1,8$, the exponent of
$e^{\tau}$ in $\ab(g)$ is equal to the exponent of $y^{\tau}$ (or
$y^{\tau,1}$, or $y^{\tau,8}$) in $g$. Therefore, the kernel
$F_{t,n}$ is generated by $\lab{1}_i$, $\lab{8}_i$ and $z$.

Recall that $n = 18$.
\begin{cor}
Let $H$ denote the group generated by $\set{z, \lab{1}_i,
\lab{8}_i}_{i = 1,\dots,n}$, with the relations \eqs{Fs0}{Fs1} and
$z$ central. Define a map $\ab \co H \ri \Z^2 = \Z e^1 \oplus \Z
e^8$ by $\ab(z) = 0$, $\ab(\lab{1}_i) = e^1$ and $\ab(\lab{8}_i) =
e^8$. The symmetric group is acting on $H$ by indices, and the
action is compatible with $\ab$.

Then $K_C = \Ker(\psi_C \co C \ri S_{n})$ is isomorphic to
$\Ker(\ab)$, and $C$ is the semidirect product $S_n \semidirect
\Ker(\ab)$ (action on the indices).
\end{cor}

Note that $H$ is an extension of $\Z^{2n} = \Z^{36}$ by $\Z
=\sg{z}$, and $\Ker(\ab)$ is an extension of $\Z^{2(n-1)} =
\Z^{34}$ by $\Z$. This proves Theorem \ref{main}.

{}Since $z$ is invariant under the action of $S_n$, it generates
the center of $S_n \semidirect H$ for the group $H$ just defined.
Therefore  $\Phi^{-1}(z)$ (or more precisely its image in $C$)
generates the center of $C$. The computations above allow us to
identify this element. Recall that $C = \sg{u_1,\dots,u_{27}}$
corresponding to the intersecting pairs of planes (with the
generators numbered as in Figure~\ref{degenerationTT}); here too
we write $j$ for $u_j$.
\begin{prop}
Let $\sigma_1 = 21 \cdot 19\cdot 8 \cdot 6$, $\tau_1 =
\sigma_1^{-1} \cdot 14 \cdot \sigma_1$, $\sigma_2 = 20 \cdot 24
\cdot 25 \cdot 16 \cdot 11 \cdot 5$, $\tau_2 = 19 \cdot 21 \cdot
14 \cdot 21 \cdot 19$, $\tau_3 = \sigma_2^{-1} \tau_2 \sigma_2$
and $\tau_4 = \sigma_2^{-1}\cdot 8\cdot \sigma_2$.

\forget $z = 6 8 19 21 14 21 19 8 6 1 5 11 16 25 24 20 19 21 14 21
19 8 20 24 25 16 11 5 4 5 11 16 25 24 20 19 21 14 21 19 20 24 25
16 11 5 1 6 8 19 21 14 21 19 8 6 5 11 16 25 24 20 19 21 14 21 19
20 24 25 16 11 5 4 5 11 16 25 24 20 8 19 21 14 21 19 20 24 25 16
11 5$ \forgotten
Then the center of $C$ is the infinite cyclic group generated by
$[\tau_1 \cdot 1 , \tau_3^{-1} (\tau_4 \cdot 4) \tau_3]$.
\end{prop}
\begin{proof}
Consider the above as elements of $\CoxY{T}$. The only generators
used which are not in $T_0$, are $1$ and $4$. Recall from \cite{v}
that $\CoxY{T_0} \isom S_{n}$ if $T_0$ is a spanning subtree. We
thus compute in the group $\CoxY{T_0} \isom S_{18}$ (numbering as
in Figure~\ref{hexT}): $\tau_1 = (2\,7)$, $\tau_2 = (7\,10)$,
$\tau_3 = (1\,7)$ and $\tau_4 = (1\,3)$. Now, $\tau_1 \cdot 1$ and
$\tau_4 \cdot 4$ are in the kernel of $\psi_T$ (see
Section~\ref{CYT} for the definition of this map), and $\Phi$ maps
them to $(\lab{1}_2)^{-1} \lab{1}_7$ and $(\lab{8}_3)^{-1}
\lab{8}_1$ respectively. Moreover $\Phi(\tau_3 (\tau_4 \cdot
4)\tau_3) = (\lab{8}_3)^{-1} \lab{8}_7$.

Now $\Phi([\tau_1 \cdot 1 , \tau_3 (\tau_4 \cdot 4) \tau_3]) =
[(\lab{1}_2)^{-1} \lab{1}_7, (\lab{8}_3)^{-1} \lab{8}_7] =
[\lab{1}_7, \lab{8}_7]$ since $\lab{8}_3$ and $\lab{1}_2$ commute
with $\lab{i}_7$ and with each other; and the last commutator is
$z$ of \Eq{Fs1}, which generates the center of $S_n \semidirect
F_{t,n}$ modulo the cyclic relations.
\end{proof}

We conclude with a general remark, motivated by a topological
interpretation of the computation done in this section.
Originally, $\CoxY{T}$ is isomorphic to $S_n$ acting on a certain
subgroup of $\pi_1(T)^n$. Adding a cyclic relation to $\CoxY{T}$
trivializes one generator, and this can be achieved by  patching a
$2$-cell (homeomorphic to $D^2$) on this cycle of $T$. The
degenerated object $X_0$ can be viewed as a triangulation of a
torus; then $T$ is the dual graph of its $1$-skeleton $S_0$.
Adding all the patches to $T$ results with a surface homeomorphic
to $X_0$, namely to the torus $\T$. The fundamental group is now
$\pi_1(\T) = \Z^2$, and indeed the  kernel of the map
$\pi_1(\T)^{n} \ri \operatorname{H}_1(\T)$ is $\Z^{2(n-1)} =
\Z^{34}$, which is the abelianization of $K_C$.

Let $X$ be a surface of general type of degree $n$, with a
degeneration to a union $X_0$ of planes where no three planes meet
in a line. In all cases computed so far (including the Hirzebruch
and Veronese surfaces, embeddings of $\C\P^1\times \C\P^1$ with
respect to the full linear system $\vert a L_1 + b L_2\vert$, as
well as $\C\P^1\times\T$ and $\T\times \T$ which is dealt with
here), the kernel $K_C = \Ker(\psi_C \co C \ri S_n)$ has the same
abelianization as the kernel of $\pi_1(X_0)^n \ri
\operatorname{H}_1(X_0)$. It would be interesting to know how far
this observation goes.

\forget A long standing problem is to construct surfaces for which
the fundamental group of the Galois cover is not virtually
solvable. Along the same lines, we suspect that this would be the
case for $X$ whenever $\pi_1(X_0)$ is not solvable. \forgotten

\forget
\section{Splitting the map $\theta$}\label{splitting}

(($\theta$ is defined in the first paragraph of the previous
section))

\begin{theorem}
The short exact sequence
$$1\llra \Ker(\theta) \llra \tilde{\pi}_1 \stackrel{\theta}{\llra} C \llra 1$$
splits.

In fact, the map $f : C \rightarrow \tilde{\pi}_1$ defined by
$f(u_1) = \Gamma_1$, ...... is a splitting.
\end{theorem}
\begin{proof}
Since $f \circ \theta$ is clearly the identity on $C$, we only
need to prove that $f$ is well defined, namely to verify the
relations described in the previous section.

Relation \eq{AX1}: % 1 \cdot 13 \cdot 22 \cdot 6 \cdot 4 \cdot 2 \cdot 4 \cdot 6  \cdot 22 \cdot 13 & = & e
\begin{eqnarray*}
f(u_1) f(u_{13}) f(u_{22}) f(u_6) f(u_4) f(u_2) f(u_4) f(u_6) f(u_{22}) f(u_{13}) & = & ... \\
    & = & ... \\
    & = & ... \\
    & = & 1.
\end{eqnarray*}
(explain each step, quoting relations from the list).

Relation \eq{AX2}: % 24 \cdot 25 \cdot 26 \cdot 25 \cdot 24 \cdot 22 \cdot 23 \cdot 27 \cdot 23 \cdot 22  & = & e\label{AX2} \\

Relation \eq{AX3}: % 5 \cdot 4 \cdot 8 \cdot 4 \cdot 5 \cdot 19 \cdot 15 \cdot 11 \cdot 15 \cdot 19  & = & e\label{AX3} \\

Relation \eq{AX4}: % 9 \cdot 10 \cdot 11 \cdot 10 \cdot 9 \cdot 16 \cdot 25 \cdot 12 \cdot 25 \cdot 16  & = & e\label{AX4} \\

Relation \eq{AX5}: % 12 \cdot 24 \cdot 20 \cdot 24 \cdot 12 \cdot 24 \cdot 20 \cdot 24 \cdot 12 \cdot 24 \cdot 20 \cdot 24  & = & e\label{AX5} \\

Relation \eq{AX6}: % 15 \cdot 21 \cdot 15 \cdot 14 \cdot 13 \cdot 14 \cdot 16 \cdot 26 \cdot 16 \cdot 14 \cdot 13 \cdot 14  & = & e\label{AX6}\\

Relation \eq{AX7}: % 19 \cdot 20 \cdot 19 \cdot 21 \cdot 19 \cdot 20 \cdot 19 \cdot 21 \cdot 19 \cdot 20 \cdot 19 \cdot 21  & = & e\label{AX7} \\

\begin{eqnarray}
20 \cdot 19 \cdot 21 \cdot 19 \cdot 20 \cdot 18 \cdot 17 \cdot 18 \cdot 27 \cdot 18 \cdot 17 \cdot 18  & = & e\label{AX8} \\
%22 \cdot 23 \cdot 22 \cdot 24 \cdot 25 \cdot 24 \cdot 22 \cdot 23 \cdot 22 \cdot 24 \cdot 25 \cdot 24  & = & e\label{AX9} \\ %-- this is a rotation of \eq{AX10}.
24 \cdot 25 \cdot 24 \cdot 22 \cdot 23 \cdot 22 \cdot 24 \cdot 25 \cdot 24 \cdot 22 \cdot 23 \cdot 22  & = & e\label{AX10} \\
3 \cdot 9 \cdot 7 \cdot 9 \cdot 3 \cdot 9 \cdot 7 \cdot 9 \cdot 3 \cdot 9 \cdot 7 \cdot 9  & = & e\label{AX11} \\
5 \cdot 2 \cdot 5 \cdot 18 \cdot 23 \cdot 18 \cdot 10 \cdot 3 \cdot 10 \cdot 18 \cdot 23 \cdot 18  & = & e\label{AX12} \\
5 \cdot 4 \cdot 5 \cdot 19 \cdot 15 \cdot 19 \cdot 5 \cdot 4 \cdot 5 \cdot 19 \cdot 15 \cdot 19  & = & e\label{AX13} \\
6 \cdot 4 \cdot 6 \cdot 22 \cdot 13 \cdot 22 \cdot 6 \cdot 4 \cdot 6 \cdot 22 \cdot 13 \cdot 22  & = & e\label{AX14}\\
6 \cdot 7 \cdot 6 \cdot 24 \cdot 20 \cdot 24 \cdot 6 \cdot 7 \cdot 6 \cdot 24 \cdot 20 \cdot 24  & = & e\label{AX15} \\
6 \cdot 7 \cdot 6 \cdot 8 \cdot 6 \cdot 7 \cdot 6 \cdot 8 \cdot 6 \cdot 7 \cdot 6 \cdot 8  & = & e\label{AX16} \\
9 \cdot 10 \cdot 9 \cdot 16 \cdot 25 \cdot 16 \cdot 9 \cdot 10 \cdot 9 \cdot 16 \cdot 25 \cdot 16  & = & e\label{AX17} \\
9 \cdot 7 \cdot 9 \cdot 14 \cdot 17 \cdot 14 \cdot 9 \cdot 7 \cdot 9 \cdot 14 \cdot 17 \cdot 14  & = & e\label{AX18} \\
1 \cdot 14 \cdot 17 \cdot 14 \cdot 9 \cdot 7 \cdot 9 \cdot 3 \cdot 9 \cdot 7 \cdot 9 \cdot 14 \cdot 17 \cdot 14  & = & e\label{AX19} \\
6 \cdot 7 \cdot 6 \cdot 8 \cdot 6 \cdot 7 \cdot 6 \cdot 24 \cdot
20 \cdot 24 \cdot 12 \cdot 24 \cdot 20 \cdot 24  & = &
e\label{AX20}
\end{eqnarray}
\begin{eqnarray}
10 \cdot 3 \cdot 10 \cdot 18 \cdot 23 \cdot 18 \cdot 10 \cdot 3
\cdot 10 \cdot 18 \cdot 23 \cdot 18 \cdot 10 \cdot \label{AX21} \\
\cdot 3 \cdot 10
\cdot 18 \cdot 23 \cdot 18  & = & e \nonumber\\
15 \cdot 14 \cdot 13 \cdot 14 \cdot 15 \cdot 21 \cdot 15 \cdot 14
\cdot 13 \cdot 14 \cdot 15 \cdot 21 \cdot 15 \cdot \label{AX22} \\
\cdot 14 \cdot 13 \cdot 14 \cdot 15 \cdot 21  & = & e \nonumber \\
19 \cdot 20 \cdot 18 \cdot 17 \cdot 18 \cdot 20 \cdot 19 \cdot 21
\cdot 19 \cdot 20 \cdot 18 \cdot 17 \cdot 18 \cdot \label{AX23}\\
\cdot 20 \cdot 19 \cdot 21 \cdot 19 \cdot 20 \cdot 18 \cdot 17
\cdot 18 \cdot 20 \cdot 19 \cdot 21  & = & e \nonumber
\end{eqnarray}
\begin{eqnarray}
25 \cdot 24 \cdot 22 \cdot 23 \cdot 22 \cdot 24 \cdot 25 \cdot 26
\cdot 25 \cdot 24 \cdot 22 \cdot 23 \cdot 22 \cdot \label{AX24}\\
\cdot 24 \cdot 25 \cdot 26 \cdot 25 \cdot 24 \cdot 22 \cdot 23
\cdot 22 \cdot 24 \cdot 25 \cdot 26  & = & e \nonumber \\
6 \cdot 4 \cdot 2 \cdot 4 \cdot 6 \cdot 22 \cdot 13 \cdot 22 \cdot
6 \cdot 4 \cdot 2 \cdot 4 \cdot 6 \cdot 22 \cdot 13 \cdot \label{AX25}\\
\cdot 22 \cdot 6 \cdot 4 \cdot 2 \cdot 4 \cdot 6 \cdot 22 \cdot 13
\cdot 22 & = & e \nonumber \\
9 \cdot 7 \cdot 9 \cdot 3 \cdot 9 \cdot 7 \cdot 9 \cdot 14
\cdot\label{AX26}
17 \cdot 14 \cdot 9 \cdot 7 \cdot 9 \cdot 3 \cdot 9 \cdot \\
7 \cdot 9 \cdot 14 \cdot 17 \cdot 14 \cdot 9 \cdot 7 \cdot 9 \cdot
3 \cdot 9 \cdot 7 \cdot 9 \cdot 14 \cdot 17 \cdot 14  & = &
e\nonumber
\end{eqnarray}

\end{proof}

\forgotten

\end{document}